\documentclass[10pt,reqno]{amsart}
\usepackage{pdfpages,color,tikz,enumitem}
\usepackage[numbers,sort&compress]{natbib}
\usepackage{graphicx}
\usepackage{amsaddr}
\usepackage{amsmath,amssymb,amsthm}
\usepackage[CJKbookmarks=True]{hyperref}

\allowdisplaybreaks[4]

\numberwithin{equation}{section}

\newtheorem{Def}{Definition}[section]
\newtheorem{Thm}[Def]{Theorem}
\newtheorem{Lem}[Def]{Lemma}

\newcommand{\R}{\mathbb{R}}
\newcommand{\J}{\mathbb{J}}
\newcommand{\mC}{\mathbb{C}}
\newcommand{\Z}{\mathbb{Z}}

\newcommand{\T}{\mathbb{T}^2}

\newcommand{\re}{\mathrm{Re}\,}
\newcommand{\im}{\mathrm{Im}\,}

\newcommand{\Hess}{\mathrm{Hess}}

\newcommand{\bvec}[1]{\boldsymbol{#1}}

\newcommand{\vj}{\bvec{j}}

\newcommand{\vq}{\bvec{q}}

\newcommand{\vx}{\bvec{x}}
\newcommand{\va}{\bvec{a}}

\newcommand{\val}{\bvec{a}^0}
\newcommand{\vb}{\bvec{b}}

\newcommand{\p}{\partial}

\newcommand{\Div}{\nabla\cdot}

\newcommand{\la}{\left\langle}
\newcommand{\ra}{\right\rangle}
\newcommand{\nn}{\nonumber\\}

\newcommand{\wto}{\rightharpoonup}

\newcommand{\be}{\begin{equation}}
\newcommand{\ee}{\end{equation}}

\newcommand{\M}{W^{-1,1}(\T)}
\newcommand{\eprint}{arxiv:}

\begin{document}
\title[Quantized Vortex Dynamics of (NLW)  on the Torus]{Quantized Vortex Dynamics of the Nonlinear Wave Equation  on the Torus}
\author[Y. Zhu]{Yongxing Zhu$^*$}
\address{Department of Mathematical Sciences, Tsinghua University, \\Beijing, 100084, China}
\thanks{This work was partially supported by the China Scholarship Council (Y. Zhu) and the National Natural Science Foundation of China [grant 12141103]. Part of the work was
done when the author was visiting National University of Singapore during 2021-2023 and the Institute for
Mathematical Science in 2023. }
\thanks{$^*$ Corresponding author: zhu-yx18@mails.tsinghua.edu.cn (Yongxing Zhu)}
\begin{abstract}
We derive rigorously the reduced dynamical laws for quantized vortex dynamics of the nonlinear wave equation on the torus when the core size of vortex $\varepsilon\to 0$. It is proved that the reduced dynamical laws are second-order nonlinear ordinary differential equations which are driven by the renormalized energy on the torus, and the initial data of the reduced dynamical laws are determined by the positions of vortices and the momentum. We will also investigate the effect of the momentum on the vortex dynamics.
\end{abstract}
\date{}
\maketitle

\section{Introduction}\label{sec:Introduction}

In this paper, we study the vortex dynamics of the nonlinear wave equation (NLW) \cite{Lin1999VortexDynamicsWave,Jerrard1999GinzburgLandauwave,BaoZengZhang2008WaveVortex}: 
\begin{equation}\label{eq:NW}
k_\varepsilon\p^2_tu^\varepsilon-\Delta u^\varepsilon+\frac{1}{\varepsilon^2}(|u^\varepsilon|^2-1)u^\varepsilon=0,\vx\in\T,t>0,
\end{equation}
with initial data 
\begin{equation}\label{eq:initial}
u^\varepsilon(\vx,0)=u_0^\varepsilon(\vx),\quad u_t^\varepsilon(\vx,0)=u^\varepsilon_1(\vx),\quad \vx\in\T.
\end{equation}
Here, $0<\varepsilon<<1$ is a parameter related to the vortex core size, $\kappa_\varepsilon=1/|\log \varepsilon|$ is a positive parameter, $\T=(\R/\Z)^2$ is the unit torus, $\vx=(x,y)^T$ is the spatial coordinate, $t$ is the time variable, $u^\varepsilon$ is a complex-valued function which is called the order parameter, and $u^\varepsilon_0,u^\varepsilon_1$ are the given complex-valued initial data. 

We define the Ginzburg-Landau functional (energy) by \cite{Lin1999VortexDynamicsWave,Jerrard1999GinzburgLandauwave}
\begin{equation}\label{eq:def of E}
E^\varepsilon(u^\varepsilon(t)):=\int_{\T}e^\varepsilon(u^\varepsilon(\vx,t))d\vx,
\end{equation}
the momentum by \cite{Jerrard1999GinzburgLandauwave}
\begin{equation}\label{eq:def of Q}
  \bvec{Q}(u^\varepsilon(t)):=\int_{\T}\vj(u^\varepsilon(\vx,t))d\vx,
\end{equation}
the Hamiltonian $h^\varepsilon(u^\varepsilon)$ by \cite{Jerrard1999GinzburgLandauwave}
\begin{equation}\label{eq:def of h}
h^\varepsilon(u^\varepsilon):=\frac{k_\varepsilon}{2}|u^\varepsilon_t|^2+e^\varepsilon(u^\varepsilon),
\end{equation}
where energy density $e^\varepsilon(v)$, current $\vj(v)$ and Jacobian $J(v)$ are defined as follows: for any complex-valued function $v:\T\to\mC^2$,
\begin{equation}\label{eq:def of e j J}
\begin{aligned}
&\vj(v):=\im (\overline{v}\nabla v), \quad
e^\varepsilon(v):=\frac{1}{2}|\nabla v|^2+\frac{1}{4\varepsilon^2}(1-|v|^2)^2,\\& J(v)=\frac{1}{2}\nabla\cdot(\J \vj(v))=\im (\partial_x\overline{v}\,\partial_y v),
\end{aligned}
\end{equation}
with $\overline{v},\im v$ denoting the complex conjugate and imaginary part of $v$, respectively, and 
\begin{equation}
  \J=\left(\begin{array}{cc}0&1\\-1&0 \end{array}\right).
\end{equation}
It's known that \eqref{eq:NW} satisfies the conservation law \cite{Jerrard1999GinzburgLandauwave}:
\begin{equation}\label{eq:hamiltonian conservation NW}
\begin{aligned}
H^\varepsilon(u^\varepsilon(t))&:=\int_{\T}h^\varepsilon(u^\varepsilon(\vx,t))d\vx=\int_{\T}\frac{k_\varepsilon}{2}|u^\varepsilon_t(\vx,t)|^2d\vx+E^\varepsilon(u^\varepsilon(t))\\&\equiv \int_{\T}\frac{k_\varepsilon}{2}|u^\varepsilon_1|^2d\vx+E^\varepsilon(u_0^\varepsilon).
\end{aligned}
\end{equation}

Quantized vortex arises in many physical phenomena, such as superfluidity, superconductivity and Bose-Einstein condensate, and is widely observed in experiments \cite{AransonKramer2002ComplexGinzburgLandau,Anderson2010ExperimentsVorticeSuperfluid}. Mathematically, quantized vortex means the topological defect of the order parameter, which is the zero point with nonzero winding number (degree) in two-dimensional cases. 
Based on mathematical analysis and numerical simulations, a quantized vortex is dynamically stable only when its degree is $\pm 1$ \cite{Mironescu1995StabilityRadialSolutionGinzburgLandau,Anderson2010ExperimentsVorticeSuperfluid,BaoZengZhang2008WaveVortex}.

The interaction and dynamics of quantized vortex were widely studied in the past decades. We refer to \cite{test,Rubinstein1991SelfinducedmotionLinedefects,Neu1990VorticesComplexScalarfield,E1994DynamicsGinzburgLandau,BaoTang2013NumericalQuantizedVortexGinzburgLandau,BaoTang2014NumericalStudyQuantizedVorticesSchrodinger,SandierSerfaty2004EstimateGinzburgLandau,Lin1996DynamicGinzburgLandau,BetheulOrlandiSmets2006GinzburgLandautoMotionByMeanCurvature,ChenSternberg2014VortexDynamicsManifold,DuJuGInzburgLandauSphereNumerical,JianLiu2006GinzburgLandauVortexMeanCurvatureFlow,Jerrard1999GinzburgLandauwave,JerrardSoner1998DynamicsGinzburgLandauVortices,JerrardSpirn2008RefinedEstimateGrossPitaevskiiVortices,SandierSerfaty2004GammaConvergenceGinzburgLandau,LinXin1999DynamicGinzburgLandauPlane,Lin1999VortexDynamicsWave,ZhangBaoDu2007SimulationVortexGinzburgLandauSchrodinger,BetheulJerrardSmets2008NLSVortexdynamicsPlane,BaoDuZhang2006RotatingBEC,LinLin2000,LinLin2003,LinLinWei2009,Neu1990} and references therein for the results on the whole plane $\R^2$ or on bounded domains with Neumman or Dirichlet boundary conditions. In particular, it was shown that for \eqref{eq:NW} on $\Omega$ (which can be taken as $\R^2$ or a bounded domain) \cite{Jerrard1999GinzburgLandauwave,Lin1999VortexDynamicsWave,BaoZengZhang2008WaveVortex} with proper $u_1^\varepsilon$ and boundary conditions, if the initial data $u^\varepsilon_0$ possesses $M$ distinct vortices $\va_1^{0\varepsilon},\cdots,\va_M^{0,\varepsilon}\in\Omega$ with degrees $d_1,\cdots,d_M\in\{\pm 1\}$ as $\varepsilon\to 0$, and $\va^{0,\varepsilon}:=(\va_1^{0\varepsilon},\cdots,$ $\va_M^{0,\varepsilon})^T$ satisfies
\begin{equation}
\lim_{\varepsilon\to 0}\va^{0,\varepsilon}= \va^0=(\va_1^0,\cdots,\va_M^0)^T,\quad \text{with}\ \va_j^0\ne \va_k^0,\quad 1\le j<k\le M,
\end{equation}
then before a collision occurs, the solution $u^\varepsilon$ also possesses $M$ distinct vortices $\va_1^\varepsilon(t),\cdots,\va_M^\varepsilon(t)$ with degrees $d_1,\cdots,d_M\in\{\pm 1\}$ as $\varepsilon\to 0$, and $\va^\varepsilon(t):=(\va_1^\varepsilon(t),\cdots,$ $\va_M^\varepsilon(t))^T$ satisfies
\begin{equation}
\lim_{\varepsilon\to 0}\va^\varepsilon(t)=\va(t)=(\va_1(t),\cdots,\va_M(t))^T,\quad \text{with}\ \va_j(t)\ne \va_M(t),\quad 1\le j<k\le M.
\end{equation}
Moreover, $\va=\va(t)$ satisfies
\[
  \ddot{\va}=-\frac{1}{\pi}\nabla W_\Omega(\va),
\]
where $W_\Omega$ is the renormalized energy on $\Omega$ of the following form:
\begin{equation}\label{eq:RDL on Omega}
  W_\Omega(\va)=-\pi\sum_{1\le j\ne k\le M}d_jd_k\log|\va_j-\va_k|+\text{Reminding term determined by} \ \Omega.
\end{equation}

However, for \eqref{eq:NW} on the torus, there are two differences: (i) The summation of degrees of vortices must be zero since $\T$ is a compact manifold \cite{CollianderJerrard1999GLvorticesSE}. Hence, we can assume $M=2N $ for an integer $N$ and
\begin{equation}
d_1=\cdots=d_N=1,\quad d_{N+1}=\cdots=d_{2N}=-1.
\end{equation}
(ii) 
The vortex dynamics on the torus are determined by both the positions of vortices and the limit momentum $\bvec{Q}_0:=\lim_{\varepsilon\to 0}\bvec{Q}(u^\varepsilon_0)$ \cite{CollianderJerrard1999GLvorticesSE,zhubaojian2023quantized,zhu2023quantized}. And if $u^\varepsilon$ satisfies
\begin{equation}\label{con:convergence of J varphi}
 J (u_0^\varepsilon)\to \pi\sum_{j=1}^{2N} d_j \delta_{\val_j} \quad \text{in}\ \M=[C^1(\T)]', \text{as}\  \varepsilon\to 0^+,
 \end{equation}
where $\delta$ is the Dirac-$\delta$ function and $\delta_{\vx_0}(\vx):=\delta(\vx-\vx_0)$ for any $\vx_0\in\T$, then \cite{zhu2023quantized}
 \begin{equation}
  \bvec{Q}_0\in 2\pi\J\sum_{j=1}^{2N}d_j\va_j^0+2\pi\Z^2.
 \end{equation}

The purpose of this article is to extend the results in \cite{Jerrard1999GinzburgLandauwave,Lin1999VortexDynamicsWave} and obtain the reduced dynamical laws of \eqref{eq:NW} involving the influence of momentum. 

Before the statement of our result, we need to introduce some notations. We define 
\[
  (\T)_*^{2N}=\{(\vx_1,\cdots,\vx_{2N})^T\in(\T)^{2N}|\vx_k\ne\vx_l \ \text{for}\ 1\le k<l\le 2N\}.
\]
For any $\va=(\va_1,\cdots,\va_{2N})^T\in(\T)_*^{2N}$ and $\vq\in 2\pi\sum_{j=1}^{2N}d_j\va_j+2\pi\Z^2$, the renormalized energy is defined by \cite{IgnatJerrard2021Renormalizedenergymanifold}
\begin{equation}\label{eq:def of W}
  W(\va;\vq):=-\pi\sum_{\begin{subarray}{c}1\le k\ne l\le 2N\end{subarray}}d_kd_lF(\va_k-\va_l)+\frac{1}{2}\left|\vq\right|^2,
\end{equation}
where $F$ is the solution of 
\begin{equation}\label{eq:define of F}
  \Delta F(\vx)=2\pi(\delta(\vx)-1),\vx\in\T,\quad\text{with}\  \int_{\T}Fd\vx=0.
\end{equation}
We define
\begin{equation}\label{eq:def of I}
    \gamma:=\lim_{\varepsilon\to 0}\left(\inf_{u\in H^1_g(B_1(\bvec{0}))}\int_{B_1(\bvec{0})}e^\varepsilon(u)d\vx-\pi\log\frac{1}{\varepsilon}\right),
\end{equation}
where $H^1_g(B_1(\bvec{0}))$ is a function space defined by
\[
  H_g^1(B_1(\bvec{0}))=\left\{u\in H^1(B_1(\bvec{0}))\left| u(\vx)=g(\vx)=\frac{x+iy}{|\vx|}\ \text{for}\  \vx\in \p B_1(\bvec{0})\right.\right\}.
\]
 Then we define 
\begin{equation}\label{eq:def of We}
  W_\varepsilon(\va):=2N\left(\pi\log\frac{1}{\varepsilon}+\gamma \right)+W(\va).
\end{equation}
Our main result is stated as follows:
\begin{Thm}[Reduced dynamical laws of the nonlinear wave equation]\label{thm:dynamics NW}
Assume that there exists $\va^0=(\va_1^0,\dots,$ $\va^0_{2N})^T\in(\T)^{2N}_*$, $\bvec{q}_0\in 2\pi\sum_{j=1}^{2N}d_j\va_j^0+2\pi\Z^2$ such that the initial data  of \eqref{eq:NW} satisfies \eqref{con:convergence of J varphi}, and
 \begin{equation}\label{con:limit of integer of j(varphi),energy}
    \bvec{Q}_0=\J\vq_0,\quad \lim_{\varepsilon\to 0}(E^\varepsilon(u^\varepsilon_0)- W_\varepsilon(\val))=0,\quad k_\varepsilon\int_{\T}|u_1^\varepsilon|^2dx\to0.
 \end{equation}
Then there exist  Lipschitz paths $\va_j:[0,T)\to\T,j=1,\cdots,2N$, such that
\[
    J(u^\varepsilon(\vx,t))\to \pi\sum_{j=1}^{2N}d_j\delta_{\va_j(t)},\quad k_\varepsilon e^\varepsilon(u^\varepsilon(\vx,t))\to\pi\sum_{j=1}^{2N}\delta_{\va_j(t)},
\]
and $\va=\va(t)=(\va_1(t),\cdots,\va_{2N}(t))^T$ satisfies
\begin{equation}\label{eq:NWODE}
  \ddot{\va}_j=-\frac{1}{\pi}\nabla_{\va_j}W(\va;\vq_*(\va)),
\end{equation}
with initial data
\begin{equation}\label{eq:initial ODE}
  \va_j(0)=\va_j^0,\quad \dot{\va}_j(0)=0,\quad \vq_*(\va(0))=\vq_0,
\end{equation}
where $\vq_*(\va)=\vq_*(\va(t))$ is continuous with respect to $t$ and satisfies
\begin{equation}\label{eq:def of q*}
  \vq_*(\va(t))\in2\pi\J\sum_{j=1}^{2N}d_j\va_j(t)+2\pi\Z^2.
  \end{equation}
\end{Thm}

The rest of this  paper is organized as follows: In Section \ref{sec:Preliminaries}, we will introduce the canonical harmonic map and the renormalized energy on the torus, and study the dynamics of the momentum of function on the torus. In Section \ref{sec:NW}, we will give the proof of Theorem \ref{thm:dynamics NW}. We begin with proving the existence of vortices and convergence of the current of $u^\varepsilon$. Then we give some lower bounds of the energy related to the solution. Finally, we prove the reduced dynamical laws. In Section \ref{sec:Numerical}, we will investigate the effect of $\bvec{Q}_0=\J\vq_0$ on the vortex dynamics.

\section{Preliminaries}\label{sec:Preliminaries}

\subsection{Canonical harmonic map and renormalized energy}

We first introduce some frequently used notations. For two complex vectors $\bvec{z}=(z_1,\cdots,z_M)^T,\bvec{w}=(w_1,\cdots,w_M)^T\in\mC^M,M=1,2$, we define
\[
    \la \bvec{z},\bvec{w}\ra:=\re\sum_{j=1}^M\overline{z}_jw_j.
\]
We use $\Hess(v)$ to denote the Hessian matrix of a function $v:\T\to\mC$.
For any $\va=(\va_1,\cdots,\va_{2N})^T\in(\T)^{2N}_*$, we define 
\begin{equation}\label{eq:def of ra}
r(\va)=\frac{1}{4}\min_{1\le j<k\le 2N}|\va_j-\va_k|.
\end{equation}
Then for any $\rho<r(\va)$, we define
\begin{equation}
\T_\rho(\va)=\T\setminus\cup_{j=1}^{2N}B_\rho(\va_j)
,\quad
\T_*(\va)=\cup_{\rho>0}\T_\rho(\va)=\T\setminus\{\va_1,\cdots,\va_{2N}\}.
\end{equation}

For any $\va=(\va_1,\cdots,\va_{2N})^T\in(\T)^{2N}_*$ and $\vq\in 2\pi\sum_{j=1}^{2N}d_j\va_j+2\pi\Z^2$, the canonical harmonic map $H=H(\vx;\va,\vq)\in C^\infty_{loc}(\T_*(\va))\cap W^{1,1}(\T)$ is defined by solving (See Lemma 2.1 and Lemma 2.2 in \cite{zhubaojian2023quantized}.)
\begin{equation}
|H|=1,\quad\Div \vj(H)=0,\quad 2J(H)=\nabla\cdot(\J \vj(H))=2\pi\sum_jd_j\delta_{\va_j},\quad \int_{\T}\vj(H)=\J\bvec{q}.\label{eq:div of jH}
\end{equation}
Then  Lemma 7 in \cite{Jerrard1999GinzburgLandauwave} implies that
for any $\eta\in C^2_0(B_{r(\va)}(\va))$ which is linear in a neighborhood of $\va_{j}$, we have
\begin{equation}
\int_{\T} \left(\la \Hess(\eta)\vj(H),\vj(H)\ra-\frac{1}{2}\Delta \eta|\vj(H)|^2\right)d\vx=-\nabla \eta(a_i)\cdot\nabla_{\va_j}W(\va;\vq).\label{eq:production of jH and eta GL}
\end{equation}

From (2.21) in \cite{zhu2023quantized}, for any small $\rho$, there exists a constant $C_\rho$ such that 
\begin{equation}\label{eq:Lip W}
\|W(\va;\vq)\|_{C^1((\T)^{2N}_\rho)}\le C_\rho,
\end{equation}
with
\begin{equation}
(\T)^{2N}_\rho:=\{(\vx_1,\cdots,\vx_{2N})\in(\T)^{2N}||\vx_j-\vx_k|>\rho,\forall 1\le j<k\le 2N\}.
\end{equation}

\subsection{Momentum of function on the torus}

\begin{Lem}\label{lem:continuity of Q}
Assume that $u^\varepsilon$ is a sequence of functions defined  over $\T\times[0,T)$ and there exists a constant $C$ independent of $t$ and $\bvec{Q}^*(t)\in L^1([0,T))$ such that
\begin{equation}\label{eq:bound of energy on slides and ut}
E^\varepsilon(u^\varepsilon(t))\le C\log\frac{1}{\varepsilon},\quad  \forall t\in [0,T),\qquad
\int_0^T\int_{\T}|\p_tu^\varepsilon|^2d\vx dt\le C\log\frac{1}{\varepsilon},
\end{equation}
and
\begin{equation}
\bvec{Q}(u^\varepsilon(t))=\int_{\T}\vj(u^\varepsilon(\vx,t))d\vx\wto \bvec{Q}^*(t)\quad \text{in}\ L^1([0,T)).
\end{equation}
 Then $\bvec{Q}^*=(Q_1^*,Q_2^*)^T:=\bvec{Q}^*(t)$ can be taken to be continuous.

 Moreover, for any $\tau\in[0,T)$ such the limiting momentum $\widetilde{\bvec{Q}}^*=\lim_{\varepsilon\to 0}\bvec{Q}(u^\varepsilon(\tau))$ exists, we have
\begin{equation}\label{eq:Q*tau=lim}
\bvec{Q}^*(\tau)=\widetilde{\bvec{Q}}^*=\lim_{\varepsilon\to 0}\bvec{Q}(u^\varepsilon(\tau)).
\end{equation}
\end{Lem}
\begin{proof}
We denote $\vj(u^\varepsilon)=(j_1(u^\varepsilon),j_2(u^\varepsilon))^T$ and define 
\begin{equation}\label{eq:def of Vk}
\alpha^\varepsilon:=\im(\overline{u^\varepsilon}\p_tu^\varepsilon),\quad \bvec{V}^\varepsilon=(V^\varepsilon_1,V^\varepsilon_2)^T:=\p_t\vj(u^\varepsilon)-\nabla\alpha^\varepsilon,
\end{equation}
and denote the set of Radon measures on $\T$ by $Ra(\T)$. 
Then \eqref{eq:bound of energy on slides and ut} and Theorem 3 in \cite{SandierSerfaty2004EstimateGinzburgLandau} imply that there exist $V_kd\vx\in L^2([0,T),Ra(\T))$ such that 
\begin{equation}\label{eq:converge of V}
\lim_{\varepsilon\to 0}\int_0^T\int_{\T}\phi V_k^\varepsilon d\vx dt=\int_0^T\int_{\T}\phi V_kd\vx dt,\quad \text{for any} \ \phi\in C^1_0(\T\times[0,T)).
\end{equation}
For any $\varphi\in C^1_0([0,T))$, substituting $\phi(\vx,t)=\varphi(t)$ and \eqref{eq:def of Vk} into \eqref{eq:converge of V}, we obtain

\vspace{-0.35cm}
\begin{align}\label{eq:derivative of Q}
&\int_0^T\int_{\T}\varphi V_kd\vx dt=\lim_{\varepsilon\to 0}\int_0^T\int_{\T}\varphi V_k^\varepsilon d\vx dt=\lim_{\varepsilon\to 0}\int_0^T\varphi\int_{\T}(\p_tj_k(u^\varepsilon)-\p_k\alpha^\varepsilon) d\vx dt\nn
&\quad=\lim_{\varepsilon\to 0}\int_0^T\varphi\p_t\left(\int_{\T}j_k(u^\varepsilon)d\vx\right) dt-\lim_{\varepsilon\to 0}\int_0^T\varphi\left(\int_{\T}\p_k\alpha^\varepsilon d\vx \right)dt\nn
&\quad=-\lim_{\varepsilon\to 0}\int_0^T\dot{\varphi}\left(\int_{\T}j_k(u^\varepsilon)d\vx\right) dt=-\int_0^T\dot{\varphi}Q_k^*dt=\int_0^T\varphi\dot{Q}_k^*dt.
\end{align}
Since $V_kd\vx\in L^2([0,T),Ra(\T))$, \eqref{eq:derivative of Q} implies that
\begin{equation}
\int_0^T\varphi\dot{Q}_k^*dt\le \|\varphi\|_{L^2([0,T))}\|V_kd\vx\|_{L^2([0,T),Ra(\T))}
\end{equation}
which together with the fact that $C^1_0([0,T))$ is dense in $L^2([0,T))$ implies $\dot{Q}_k^*\in L^2([0,T))$. Hence, $\bvec{Q}^*=(Q_1^*,Q_2^*)^T\in H^1([0,T))\hookrightarrow\hookrightarrow C([0,T))$ and $\bvec{Q}^*$ is continuous.

For any $\tau\in[0,T)$ such that the limiting momentum exists, we consider $U^\varepsilon\in H^1(\T\times[\tau-1,T))$ defined by
\begin{equation}
U^\varepsilon(\vx,t)=\left\{\begin{array}{ll}
u^\varepsilon(\vx,t),&t\ge \tau,\\
u^\varepsilon(\vx,\tau),&t<\tau.
\end{array}\right.
\end{equation} 
Clearly, $U^\varepsilon$ also satisfies \eqref{eq:bound of energy on slides and ut} and 
\begin{equation}\label{eq:Q*tau}
\bvec{Q}(U^\varepsilon(t))\wto \widetilde{\bvec{Q}}^*(t):=\left\{\begin{array}{ll}
\bvec{Q}^*(t),&t\ge \tau,\\
\widetilde{\bvec{Q}}^*,&t<\tau.\end{array}\right.
\end{equation}
Hence, the proof of the continuity of $\bvec{Q}^*$ also works for $\widetilde{\bvec{Q}}^*$. Combining continuity of $\widetilde{\bvec{Q}}^*$ and \eqref{eq:Q*tau}, we obtain \eqref{eq:Q*tau=lim}.
\end{proof}

\section{Vortex Dynamics}\label{sec:NW}
We first give some useful equalities related to \eqref{eq:NW}.
Any solution $u^\varepsilon$ of equation \eqref{eq:NW} satisfies the following equalities \cite{Lin1999VortexDynamicsWave,Jerrard1999GinzburgLandauwave}:

\vspace{-0.35cm}
\begin{align}
\frac{\p}{\p t}h^\varepsilon(u^\varepsilon(\vx,t))&=\nabla\cdot(\re(\overline{u^\varepsilon_t}\nabla u^\varepsilon)),\label{eq:derivative of h NW}\\
\nabla\cdot \vj(u^\varepsilon(\vx,t))&=k_\varepsilon\p_t\im(\overline{u^\varepsilon}u^\varepsilon_t),\label{eq:derivative of j NW}
\end{align}
For any $\varphi\in C^\infty(\T)$, we have \cite{Jerrard1999GinzburgLandauwave}
\begin{equation}\label{eq:mainequalityNW}
    k_\varepsilon\int_{\T}\varphi\p^2_{t}h^\varepsilon(u^\varepsilon)d\vx=\int_{\T}\left(\la\Hess(\varphi)\nabla u^\varepsilon,\nabla u^\varepsilon\ra+\left(\frac{k_\varepsilon}{2}|u_t^\varepsilon|^2-e^\varepsilon(u^\varepsilon)\right)\Delta\varphi \right)d\vx.
\end{equation}

\subsection{Convergences of the Jacobian and the current of the solution}

\begin{Lem}\label{thm: existence of vortices NW}
If the initial data $u^\varepsilon_0$ of \eqref{eq:NW} satisfies  \eqref{con:convergence of J varphi} and \eqref{con:limit of integer of j(varphi),energy}, then there exist $C^{1,1}$-paths $\vb_j:[0,T)\to\T, j=1,\cdots,2N$, such that
\begin{equation}\label{eq:converge of Ju and eu NW}
    J(u^\varepsilon(\vx,t))\to \pi\sum_{j=1}^{2N}d_j\delta_{\vb_j(t)},\quad k_\varepsilon e^\varepsilon(u^\varepsilon(\vx,t))\to\pi\sum_{j=1}^{2N}\delta_{\vb_j(t)}\quad \text{in}\ \M.
\end{equation}
\end{Lem}

\begin{proof}
We denote
\begin{equation}\label{eq:def of r0}
    r_0=r(\va^0)=\frac{1}{4}\min_{j\ne k}|\va_j^0-\va_k^0|>0.
\end{equation}
For $\varepsilon$ small, assumption \eqref{con:convergence of J varphi} implies that
\begin{equation}\label{eq:diff bet Ju0 and dirac NW}
\left\|J(u_0^\varepsilon)- \pi\sum_{j=1}^{2N} d_j \delta _{\va_{j}} \right\|_{\M}\le \frac{\pi}{400}r_0.
\end{equation}
For each small $\varepsilon$, we know $u^\varepsilon(t)\to u_0^\varepsilon$ as $t\to 0$ in $H^1(\T)$, so 
\[J(u^\varepsilon(t))\to J (u_0^\varepsilon)\quad \text{in}\ L^1(\T)\ \text{and}\ \M.
\]
 We define
\begin{equation}\label{eq:define of Te NW}
T^\varepsilon:=\sup\left\{T>0:\|J(u_0^\varepsilon)- J(u^\varepsilon(t)) \|_{\M}\le \frac{\pi}{400}r_0\ \text{for any}\ t\in[0,T]\right\},
\end{equation}
which together with \eqref{eq:diff bet Ju0 and dirac NW} implies
\begin{equation}\label{eq:diff bet Jut and dirac NW}
\left\|J(u^\varepsilon(t))- \pi\sum_{j=1}^{2N} d_j \delta _{\va_{j}} \right\|_{\M}\le \frac{\pi}{200}r_0\quad \text{for any}\ t<T^\varepsilon.
\end{equation}
\eqref{con:limit of integer of j(varphi),energy}, \eqref{eq:hamiltonian conservation NW} and \eqref{eq:def of We} imply that there exists a constant $C>0$ such that for all $\varepsilon>0$,
\begin{equation}\label{con:weak bd of energy NW}
E^\varepsilon(u^\varepsilon(t))\le H^\varepsilon(u^\varepsilon(t))\le  2N \pi\log\left(\frac{1}{\varepsilon} \right)+C.
 \end{equation}
By Theorem 1.4.4 in \cite{CollianderJerrard1999GLvorticesSE}, \eqref{eq:diff bet Jut and dirac NW} and \eqref{con:weak bd of energy NW} imply that for each $t>0,\varepsilon<\varepsilon_0$, there exist $\vb_j^\varepsilon(t)\in\T,j=1,\cdots,2N$,  such that

\vspace{-0.35cm}
\begin{align}
&\left\|J(u^\varepsilon(\vx,t))-\pi\sum_{j=1}^{2N}d_j\delta_{\vb_j^\varepsilon(t)} \right\|_{\M}=o(1),\label{eq:distance bet Jut and dirac}\\
&\left\|k_\varepsilon e^\varepsilon(u^\varepsilon(\vx,t))-\pi\sum_{j=1}^{2N} \delta_{\vb_j^\varepsilon(t)}\right\|_{\M}=o(1),\label{eq:distance bet eut and dirac}
\end{align}

\vspace{-0.30cm}
\begin{align}
\int_{\T_{r_0}(\va^0)}e^\varepsilon(u^\varepsilon(\vx,t))d\vx \le C,\quad \|\vj(u^\varepsilon)\|_{L^1(\T)}\le C,\label{eq:bound of ju, eu out balls NW}
\end{align}

\vspace{-0.30cm}
\begin{align}
d E^\varepsilon(u^\varepsilon(t))\ge 2N\pi\log\frac{1}{\varepsilon}-C.\label{eq:lower bound of energy t NW}
\end{align}
\eqref{eq:lower bound of energy t NW} and \eqref{con:weak bd of energy NW} imply
\begin{equation}\label{eq:bound of ut NW}
    k_\varepsilon \|u^\varepsilon_t\|_{L^2(\T)}^2\le C
\end{equation}
and 
\begin{equation}\label{eq:distance bet eut and hut NW}
    \|k_\varepsilon e^\varepsilon(u^\varepsilon(\vx,t))-k_\varepsilon h^\varepsilon(u^\varepsilon(\vx,t)) \|_{\M}=k_\varepsilon^2\left\||u_t^\varepsilon(\vx,t)|^2\right\|_{\M}\le Ck_\varepsilon.
\end{equation}

Then we estimate $|\vb_j^\varepsilon(t_2)-\vb_j^\varepsilon(t_1)|$: We find $\eta\in C^\infty_0(B_{r_0}(\va^\varepsilon_j(t)))$ such that
\[
    \eta(\vx)=(\vx-\vb_j^\varepsilon(t_1))\cdot\frac{\vb_j^\varepsilon(t_2)-\vb_j^\varepsilon(t_1)}{|\vb_j^\varepsilon(t_2)-\vb_j^\varepsilon(t_1)|},\quad \vx\in B_{r_0}(\vb_j^\varepsilon(t_1))
\]
then combining \eqref{eq:distance bet eut and dirac}, \eqref{eq:distance bet eut and hut NW}, \eqref{eq:derivative of h NW}, \eqref{con:weak bd of energy NW} and \eqref{eq:bound of ut NW}, we obtain

\vspace{-0.35cm}
\begin{align*}
&|\vb_j^\varepsilon(t_2)-\vb_j^\varepsilon(t_1)|=\int_{t_1}^{t_2}\int k_\varepsilon\eta \p_th^\varepsilon d\vx dt+o(1)\\&\quad=-\int_{t_1}^{t_2}\int_{\T}k_\varepsilon\la\nabla \eta,\re(\overline{u_t^\varepsilon}\nabla u^\varepsilon)\ra d\vx dt+o(1)\\
&\quad\le\|\nabla \eta\|_{L^\infty(\T)}k_\varepsilon\int_{t_1}^{t_2}\|\nabla u^\varepsilon\|_{L^2(\T)}\|u^\varepsilon_t\|_{L^2(\T)}dt+o(1)\le C|t_2-t_1|+o(1).
\end{align*}
Hence, we can find Lipschitz paths $\vb_j$'s and $T_0>0$ such that $\vb_j^\varepsilon(t)\to \vb_j(t)$ in $[0,T_0]$, which together with \eqref{eq:distance bet Jut and dirac} and \eqref{eq:distance bet eut and dirac}  implies \eqref{eq:converge of Ju and eu NW}. In particular, \eqref{con:convergence of J varphi} and \eqref{eq:converge of Ju and eu NW} imply
\begin{equation}\label{eq:bj0=aj0}
    \vb_j(0)=\va_j^0.
\end{equation}

Then, we prove that $\vb_j\in C^{1,1}([0,T_0],\T)$.
We define
\begin{equation}\label{eq:def of zeta h}
    \zeta^h(t)=\begin{cases}
    0,&|t|\ge h,\\
    \frac{1}{h}-{|t|}{h},&|t|\le h.
    \end{cases}
\end{equation}
And for any unit vector $\bvec{\nu}$, we find $\eta\in C^\infty_0(B_{r_0}(\vb_j(t)))$ satisfying 
\[
    \eta(\vx)=(\vx-\vb_j(t))\cdot\bvec{\nu},\quad \vx\in B_{3r_0/4}(\vb_j(t)).
\]
Then \eqref{eq:converge of Ju and eu NW}, \eqref{eq:distance bet eut and hut NW}, \eqref{eq:def of zeta h}, \eqref{eq:mainequalityNW}, \eqref{eq:bound of ju, eu out balls NW}, and \eqref{eq:bound of ut NW} imply

\vspace{-0.35cm}
\begin{align}\label{eq:integration of h eta}
&\frac{\pi}{h^2}(\vb_j(t-h)-2\vb_j(t)+\vb_j(t+h))\cdot\bvec{\nu}=\frac{\pi}{h^2}(\eta(\vb_j(t-h))-2\eta(\vb_j(t))+\eta(\vb_j(t+h)))\nn&\quad=\lim_{\varepsilon\to0}\frac{1}{h^2}\int_{\T}\eta k_\varepsilon(h^\varepsilon(u^\varepsilon(\vx,t-h))-2h^\varepsilon(u^\varepsilon(\vx,t))+h^\varepsilon(u^\varepsilon(\vx,t+h)))d\vx\nn
&\quad=\lim_{\varepsilon\to0}\int_{\R}\zeta^h(t-s)\int_{\T}\eta k_\varepsilon\p_t^2h^\varepsilon(u^\varepsilon(\vx,s))d\vx ds\nn
&\quad=\lim_{\varepsilon\to0}\int_{\R}\zeta^h(t-s)\int_{\T}\left(\la\Hess(\eta)\nabla u^\varepsilon,\nabla u^\varepsilon \ra +\Delta\eta\left(\frac{k_\varepsilon}{2}|u_t^\varepsilon|^2-e^\varepsilon(u^\varepsilon) \right)\right)d\vx ds\nn
&\quad\le\|\eta\|_{C^2(\T)}\limsup_{\varepsilon\to0}C\sup_{s\in[t-h,t+h]}\int_{B_{r_0}\setminus B_{3r_0/4}(\vb_j(t))}h^\varepsilon(u^\varepsilon(\vx,s))d\vx \le C.
\end{align}
Substituting $\bvec{\nu}=\frac{\vb_j(t-h)-2\vb_j(t)+\vb_j(t+h)}{|\vb_j(t-h)-2\vb_j(t)+\vb_j(t+h)|}$ (or any unit vector if $\vb_j(t-h)-2\vb_j(t)+\vb_j(t+h)=\bvec{0}$) into \eqref{eq:integration of h eta}, we have that for any $t\in[0,T_0]$ and $0<h<<1$, 
\[
    \frac{\pi}{h^2}|\vb_j(t-h)-2\vb_j(t)+\vb_j(t+h)|\le C,
\]
which implies that $\vb_j\in C^{1,1}([0,T_0],\T)$.

The proof above can be repeated as long as $\vb_j(t)\ne \vb_k(t)$ for any $1\le j<k\le 2N$. So we can extend the lifespan of $\vb_j$, i.e. $[0,T_0]$ to $[0,T)$ with 
\[
    T:=\inf\{t>0|\vb_j(t)=\vb_k(t)\ \text{for some}\ 1\le j<k\le 2N\}.\qedhere
\]
\end{proof}

With $\vb=\vb(t)=(\vb_1(t),\cdots,\vb_{2N}(t))^T$ obtained in Lemma \ref{thm: existence of vortices NW}, we can define
\begin{equation}\label{eq:def of u* NW}
    u^*:=u^*(\vx,t)=H(\vx;\vb(t),\vq_*(\vb(t))),
\end{equation}
where $\vq_*(\vb(t))$ is continuous and satisfies
\begin{equation}\label{eq:def of qb}
    \vq_*(\vb(0))=\vq_0,\quad \vq_*(\vb(t))\in2\pi\sum_{j=1}^{2N}d_j\vb_j(t)+2\pi\Z^2,
\end{equation}
and
 $H(\vx;\vb(t),\vq_*(\vb(t)))$ is defined by \eqref{eq:div of jH}.

\begin{Lem}\label{lem:converge of j j/|| NW}
Assume that $u^\varepsilon, \vb$ are the same as in Lemma \ref{thm: existence of vortices NW}. Then for any $T_1<T$ and $\rho<r_1:=\min_{t\in[0,T_1]}r(\vb(t))$, we have
\begin{equation}\label{eq:def of j and j/|| NW}
\vj(u^\varepsilon)\wto\vj(u^*)\ \text{in}\ L^1(\T\times[0,T_1]),\quad \frac{\vj(u^\varepsilon)}{|u^\varepsilon|}\wto \vj(u^*) \ \text{in}\ L^2(\T_\rho(\vb(t))\times[0,T_1]).
\end{equation}
\end{Lem}

\begin{proof}
\eqref{eq:bound of ju, eu out balls NW} implies that
 there exists a function $\vj_1\in L^1(\T\times[0,T_1])$ such that
\begin{equation}\label{eq:converge of j NW}
\vj(u^\varepsilon)\wto \vj_1\quad \text{in}\ L^1(\T\times[0,T_1]).
\end{equation}
Moreover, Theorem 1.4.4 in \cite{CollianderJerrard1999GLvorticesSE}, \eqref{con:weak bd of energy NW} and \eqref{eq:converge of Ju and eu NW} imply that for $\varepsilon$ small enough, we have
\begin{equation}\label{eq:energy out balls radius rho NW}
    \int_{\T_\rho(\vb(t))}e^\varepsilon(u^\varepsilon(\vx,t))d\vx\le C.
\end{equation}
Combining \eqref{eq:energy out balls radius rho NW} and \eqref{eq:def of e j J}, we have

\vspace{-0.35cm}
\begin{align}
\int_{0}^{T_1}\left\|\frac{\vj(u^\varepsilon)}{|u^\varepsilon|} \right\|_{L^2(\T_\rho(\vb(t)))}^2dt\le \int_{0}^{T_1}\int_{\T_\rho(\vb(t))}e^\varepsilon(u^\varepsilon(t))d\vx dt\le CT_1,\label{eq:bound of j/|| out balls NW}
\end{align}
\begin{align}
\left\||u^\varepsilon|^2-1\right\|_{L^2(\T_{\rho}(\vb(t)))}\le C\varepsilon.\label{eq:bound of L2norm NW}
\end{align}
\eqref{eq:bound of j/|| out balls NW} implies that there exists a function $\vj_2\in L^2(\T_\rho(\vb(t))\times[0,T_1])$ such that 
\begin{equation}\label{eq:converge of j/|| to j2}
\frac{\vj(u^\varepsilon)}{|u^\varepsilon|}\wto \vj_2\quad \text{in}\ L^2(\T_\rho(\vb(t))\times[0,T_1]),
\end{equation}
which together with \eqref{eq:bound of L2norm NW} implies that 
\begin{equation}\label{eq:weak converge of j out balls NW}
\vj(u^\varepsilon)=|u^\varepsilon|\frac{\vj(u^\varepsilon)}{|u^\varepsilon|}\wto \vj_2\quad \text{in}\ L^1(\T_\rho(\vb(t))\times[0,T_1]).
\end{equation}
Combining \eqref{eq:converge of j NW} and \eqref{eq:weak converge of j out balls NW}, we have 
\begin{equation}\label{eq:j1=j2}
    \vj_1=\vj_2.
\end{equation}

For any $\varphi\in C^\infty_0(\T\times[0,T_1])$, combing \eqref{eq:converge of j NW}, \eqref{eq:derivative of j NW}, \eqref{eq:bound of L2norm NW} and \eqref{eq:bound of ut NW}, we have

\vspace{-0.35cm}
\begin{align*}
    &\left|\int_{\T\times[0,T_1]}\nabla\varphi\cdot\vj d\vx dt\right|=\lim_{\varepsilon\to0}\left|\int_{\T\times[0,T_1]}\nabla\varphi\cdot\vj(u^\varepsilon)d\vx dt\right|\\&\quad=\lim_{\varepsilon\to0}\left|\int_{\T\times[0,T_1]}\varphi\nabla\cdot\vj(u^\varepsilon)d\vx dt\right|
    =\lim_{\varepsilon\to 0}\left|\int_{\T\times[0,T_1]} k_\varepsilon \varphi\frac{\p}{\p t}\im(\overline{u^\varepsilon}u^\varepsilon_t)d\vx dt\right|\\&\quad=\lim_{\varepsilon\to 0}\left|\int_{\T\times[0,T_1]} k_\varepsilon \frac{\p}{\p t}\varphi\im(\overline{u^\varepsilon}u^\varepsilon_t)d\vx dt\right|\le\lim_{\varepsilon\to 0} Ck_\varepsilon\|u_t^\varepsilon\|_{L^2(\T\times[0,T_1])} =0,
\end{align*}
which implies $\nabla\cdot \vj_1=0$. Then \eqref{eq:div of jH} and \eqref{eq:def of u* NW} imply
\[
    \nabla \cdot (\vj_1-\vj(u^*))=0.
\]
Similarly, we have $\nabla\cdot(\J(\vj_1-\vj(u^*)))=0$.
Hence, we can find $\bvec{g}:[0,T_1]\to\R^2$ such that
\begin{equation}\label{eq:def of g}\bvec{g}(t)=\vj_1(\vx,t)-\vj(u^*(\vx,t)).\end{equation} 
Noting \eqref{eq:converge of j NW}, \eqref{eq:div of jH} and \eqref{eq:def of u* NW}, we have
\begin{equation}\label{eq:intgeration of j1-ju*}
\begin{aligned}
&\bvec{g}(t)=\int_{\T}(\vj_1(\vx,t)-\vj(u^*(\vx,t)))d\vx=\int_{\T}\vj_1(\vx,t)d\vx-\J\vq_*(\vb(t)) ,\\ &\bvec{Q}(u^\varepsilon(t))\wto\int_{\T}\vj_1(\vx,t)d\vx=\bvec{g}(t)+\J\vq_*(\vb(t))\quad \text{in}\ L^1([0,T_1]).
\end{aligned}
\end{equation}
\eqref{con:weak bd of energy NW}, \eqref{eq:bound of ut NW}, \eqref{eq:intgeration of j1-ju*}, and Lemma \ref{lem:continuity of Q} imply that $\bvec{g}(t)$ is continuous with respect to $t$, and that for any $t$ such that $\lim_{\varepsilon\to 0}\bvec{Q}(u^\varepsilon(t))$ exists, we have
\begin{equation}\label{eq:gt=lim NW}
\bvec{g}(t)=\lim_{\varepsilon\to 0}\bvec{Q}(u^\varepsilon(t))-\J\vq_*(\vb(t)).
\end{equation}
Noting that for any $t\in[0,T_1]$ and any subsequence of $\bvec{Q}(u^\varepsilon(t))$ which is denoted by $\bvec{Q}(u^{\varepsilon_n}(t))$, \eqref{eq:bound of ju, eu out balls NW} and \eqref{eq:gt=lim NW} imply that we can extract a subsequence of $\bvec{Q}(u^{\varepsilon_n}(t))$ such that $\lim_{\varepsilon_n\to 0}\bvec{Q}(u^{\varepsilon_n}(t))$ exists and is equal to $\bvec{g}(t)+\J\vq_*(\vb(t))$, which implies that \eqref{eq:gt=lim NW} holds for any $t\in[0,T_1]$.
Combining \eqref{con:limit of integer of j(varphi),energy}, \eqref{con:weak bd of energy NW}, \eqref{eq:converge of Ju and eu NW}, Lemma 2.2 in \cite{zhu2023quantized}, we have
\begin{equation}
\lim_{\varepsilon\to 0}\bvec{Q}(u^\varepsilon(0))=\J\vq_0,\quad \lim_{\varepsilon\to 0}\bvec{Q}(u^\varepsilon(t))\in 2\pi\J\sum_{j=1}^{2N}d_j\vb_j(t)+2\pi\Z^2,
\end{equation}
which together with \eqref{eq:gt=lim NW} and \eqref{eq:def of qb} implies
\begin{equation}\label{eq:gt in Z2}
    \bvec{g}(0)=\bvec{0},\quad \bvec{g}(t)\in 2\pi\Z^2.
\end{equation}
Noting that $\bvec{g}$ is continuous, \eqref{eq:gt in Z2} gives $\bvec{g}(t)\equiv \bvec{0}$, which together with \eqref{eq:def of g}, \eqref{eq:converge of j NW}, \eqref{eq:converge of j/|| to j2} and \eqref{eq:j1=j2} implies \eqref{eq:def of j and j/|| NW}.
\end{proof}

\subsection[Lower bounds of E(u) and the L2 norm of ut]{Lower bounds of $E^\varepsilon(u^\varepsilon(t))$ and the $L^2$-norm of $u^\varepsilon_t$ }
\eqref{eq:converge of Ju and eu NW} and \eqref{eq:def of j and j/|| NW} together with Lemma 2.3 in \cite{zhu2023quantized} imply that for any $0\le t_1<t_2<T$, and $\rho<\min_{t\in[t_1,t_2]}r(\vb(t))$, we have
\begin{equation}\label{eq:distance between ju and ju*}
\limsup_{\varepsilon\to 0}\int_{t_1}^{t_2}\int_{\T_\rho(\vb(t))}\left(e^\varepsilon(|u^\varepsilon|)+\frac{1}{2}\left|\frac{\vj(u^\varepsilon)}{|u^\varepsilon|}-\vj(u^*) \right| \right)d\vx dt\le C\int_{t_1}^{t_2}\Sigma(t)dt,
\end{equation}
with
\[
    \Sigma(t):=\limsup_{\varepsilon\to 0}\left(E^\varepsilon(u^\varepsilon(t))-W_\varepsilon(\vb(t);\vq_*(\vb(t))) \right).
\]
Since the left hand-side of \eqref{eq:distance between ju and ju*} is always nonnegative, we have $\Sigma(t)\ge 0$, i.e.
\begin{equation}\label{eq:refined lower bound of e}
\limsup_{\varepsilon\to 0}(E^\varepsilon(u^\varepsilon(t))- W_\varepsilon(\vb(t);\vq_*(\vb(t))))\ge 0.
\end{equation}

We denote $\Omega_t=\cup_{j=1}^{2N}B_\rho(\vb_j(t))$. Then we can take $|t_2-t_1|$  small enough such that $\va_j(s)\in \Omega_t$ for any $s,t$$\in[t_1,t_2]$, which together with \eqref{eq:converge of Ju and eu NW} implies
\begin{equation}\label{eq:converge of Ju Omegat}
    J(u^\varepsilon(s))\wto \pi\sum_{j=1}^{2N}d_j\delta_{\vb_j(s)}\quad\text{in}\ W^{-1,1}(\Omega_t).
\end{equation}
Combining $\Omega_t\subset\T$, \eqref{eq:bound of ut NW} and \eqref{con:weak bd of energy NW}, we have that for any $s,t\in[t_1,t_2]$,
\begin{equation}\label{eq:energy bound on Omegat}
    \int_{\Omega_t\times[t_1,t_2]}|u_t^\varepsilon|^2d\vx ds\le C|\log \varepsilon|,\quad \int_{\Omega_t}e^\varepsilon(u^\varepsilon(s))d\vx\le 2N\pi|\log \varepsilon|+C.
\end{equation}
Then Corollary 7 in \cite{SandierSerfaty2004EstimateGinzburgLandau}, \eqref{eq:converge of Ju Omegat} and \eqref{eq:energy bound on Omegat} imply 

\vspace{-0.35cm}
\begin{align}\label{eq:lower bound of L2 ut}
    &\liminf_{\varepsilon\to 0}k_\varepsilon\int_{\T\times[t_1,t_2]}|u^\varepsilon_t(\vx,s)|^2d\vx ds\nn&\quad\ge \liminf_{\varepsilon\to 0}k_\varepsilon\int_{\Omega_t\times[t_1,t_2]}|u^\varepsilon_t(\vx,s)|^2d\vx ds\ge \pi\sum_{j=1}^{2N}\int_{t_1}^{t_2}|\dot{\vb}_j(s)|^2ds.
\end{align}

\subsection{Proof of Theorem \ref{thm:dynamics NW}}

\begin{proof}
We assume that $\va$ is the solution of \eqref{eq:NWODE}, and recall that $\vb$ is obtained in Lemma \ref{thm: existence of vortices NW}. Then we define
\begin{equation}\label{eq:def of zeta}
    \zeta(t)=\sum_{j=1}^{2N}(|\vb_j(t)-\va_j(t)|+|\dot{\vb}_j(t)-\dot{\va}_j(t)|).
\end{equation}
And we can find $T_*<T$ such that 
\begin{equation}\label{eq:def of r2 and T2}
  \sup_{t\in[0,T_*]}\zeta(t)<r_*:=\inf_{t\in[0,T_*]}\min\{r(\va(t)),r(\vb(t))\}.
\end{equation}
For simplicity, we will use the following notations in this subsection: 
\begin{equation}\label{eq:simple notation}
\begin{array}{ll}
    W(\va(t))=W(\va(t),\vq_*(\va(t))),&\quad W(\vb(t))=W(\vb(t),\vq_*(\vb(t))),\\
    W_\varepsilon(\va(t))=W_\varepsilon(\va(t),\vq_*(\va(t))),&\quad W_\varepsilon(\vb(t))=W_\varepsilon(\vb(t),\vq_*(\vb(t))).
\end{array}
\end{equation}

We first verify that $\dot{\vb}(0)=\bvec{0}$. 
\eqref{eq:lower bound of L2 ut},  \eqref{eq:hamiltonian conservation NW}, \eqref{con:limit of integer of j(varphi),energy}, \eqref{eq:refined lower bound of e} and \eqref{eq:bj0=aj0} imply

\vspace{-0.35cm}
\begin{align}\label{eq:vbdot=0}
&\frac{\pi}{2h}\int_0^h|\dot{\vb}(t)|^2dt\le\liminf_{\varepsilon\to0}\frac{1}{h}\int_0^h\int_{\T}\frac{k_\varepsilon}{2}|u_t^\varepsilon|^2 d\vx dt\nn&\quad=\liminf_{\varepsilon\to0}\frac{1}{h}\int_0^h\int_{\T}(h^\varepsilon(u^\varepsilon(t))-e^\varepsilon(u^\varepsilon(t)))d\vx dt\nn
&\quad\le\frac{1}{h}\int_0^h\left(W_\varepsilon(\va^0)-W_\varepsilon(\vb(t)) \right)dt\le C\frac{1}{h}\int_0^h|\va_j^0-\vb_j(t)|dt\le Ch.
\end{align}
Letting $h\to 0$ on the both sides of \eqref{eq:vbdot=0}, we obtain $\dot{\vb}(0)=\bvec{0}$, which together with \eqref{eq:def of zeta}, \eqref{eq:initial ODE} and \eqref{eq:bj0=aj0} implies
\begin{equation}\label{eq:zeta0=0}
    \zeta(0)=0.
\end{equation}

Taking the derivative of $\zeta$ with respect to $t$ for $t\in[0,T_*]$, and noting \eqref{eq:NWODE}, \eqref{eq:def of zeta} and \eqref{eq:Lip W}, we have

\vspace{-0.35cm}
\begin{align}\label{eq:first estimate of zeta' NW}
    \dot{\zeta}(t)\le&\sum_{j=1}^{2N}|\dot{\vb}_j(t)-\dot{\va}_j(t)|+|\ddot{\vb}_j(t)-\ddot{\va}_j(t)|\nn
    \le&\zeta(t)+\frac{1}{\pi}\sum_{j=1}^{2N}\left|\nabla_{\vb_j}W(\vb(t))+\pi\ddot{\vb}_j(t) \right|+\frac{1}{\pi}\sum_{j=1}^{2N}\left|\nabla_{\va_j}W(\va(t))-\nabla_{\vb_j}W(\vb(t)) \right|\nn
    \le& \frac{1}{\pi}\sum_{j=1}^{2N}|\bvec{A}_j(t)|+C\zeta(t).
\end{align}
where
\begin{equation}\label{eq:def of Aj}
    \bvec{A}_j(t)=\nabla_{\vb_j}W(\vb(t))+\pi\ddot{\vb}_j(t).
\end{equation}
We take $\bvec{\nu}=\bvec{A}_j(t)/|\bvec{A}_j(t)|$ (or any unit vector if $\bvec{A}_j(t)=\bvec{0}$), which together with \eqref{eq:def of Aj} gives
\begin{equation}\label{eq:Aj dot nu}
     \bvec{A}_j(t)\cdot\bvec{\nu}=|\bvec{A}_j(t)|.
 \end{equation} Then we can find $\eta\in C^\infty_0(B_{r_*}(\vb_j(t)))$ satisfying
\begin{equation}\label{eq:def and gradient of eta}
    \eta(\vx)=(\vx-\vb_j(t))\cdot\bvec{\nu},\quad  \nabla\eta(\vx)=\bvec{\nu},\quad \vx\in B_{3r_*/4}(\vb_j(t)).
\end{equation}
Recall that $\zeta^h$ is defined by \eqref{eq:def of zeta h}. 
Then \eqref{eq:integration of h eta} and \eqref{eq:def and gradient of eta} imply

\vspace{-0.35cm}
\begin{align}\label{eq:nu dot vbddot}
&\pi\ddot{\vb}_j(t)\cdot\bvec{\nu}=\lim_{h\to 0}\frac{\pi}{h^2}(\vb_j(t+h)-2\vb_j(t)+\vb_j(t-h))\cdot\bvec{\nu}\nn 
&\quad=\lim_{h\to 0}\lim_{\varepsilon\to 0}\int_{\R}\zeta^h(t-s) \int_{\T}\left(\la\Hess(\eta)\nabla u^\varepsilon,\nabla u^\varepsilon \ra+\Delta\eta\left(\frac{k_\varepsilon}{2}|u_t^\varepsilon|^2-e^\varepsilon(u^\varepsilon) \right) \right)d\vx ds.
\end{align}
\eqref{eq:production of jH and eta GL}, \eqref{eq:def of u* NW} and \eqref{eq:def and gradient of eta} imply 

\vspace{-0.35cm}
\begin{align}\label{eq:nu dot nabla W}
    &\bvec{\nu}\cdot\nabla_{\vb_j}W(\vb(t))=\lim_{h\to0}\int_{\R}\zeta^h(t-s)\nabla\eta(\vb_j(s))\cdot\nabla_{\vb_j}W(\vb_j(s))ds\nn 
    &\quad=\lim_{h\to0}\int_{\R} \zeta^h(t-s)\int_{\T}\left(\la\Hess(\eta)\vj(u^*),\vj(u^*) \ra-\Delta\eta\frac{1}{2}|\vj(u^*)|^2\right)d\vx ds.
\end{align}
Combining \eqref{eq:def of Aj}, \eqref{eq:Aj dot nu}, \eqref{eq:nu dot vbddot}, \eqref{eq:nu dot nabla W}, and using (4.17) in \cite{zhu2023quantized}, we have

\vspace{-0.35cm}
\begin{align}\label{eq:splitting of A NW}
&\left|\bvec{A}_j(t) \right|=\pi\ddot{\vb}_j(t)\cdot\bvec{\nu}+\bvec{\nu}\cdot\nabla_{\vb_j}W(\vb(t))\\
&\quad=\lim_{h\to 0}\lim_{\varepsilon\to 0}\int_{\R}\zeta^h(t-s) \int_{\T}\left(\la\Hess(\eta)\nabla u^\varepsilon,\nabla u^\varepsilon \ra+\Delta\eta\left(\frac{k_\varepsilon}{2}|u_t^\varepsilon|^2-e^\varepsilon(u^\varepsilon) \right) \right)d\vx ds\nn
&\qquad-\lim_{h\to0}\int_{\R} \zeta^h(t-s)\int_{\T}\left(\la\Hess(\eta)\vj(u^*),\vj(u^*) \ra-\Delta\eta\frac{1}{2}|\vj(u^*)|^2\right)d\vx ds\nn
&\quad=I_1+I_2+I_3
\end{align}
with 

\vspace{-0.35cm}
\begin{align}\label{eq:def of I1}
I_1=&\lim_{h\to0}\lim_{\varepsilon\to 0}\int_{\R} \zeta^h(t-s)\int_{\T}\left(\la\Hess(\eta)\nabla|u^\varepsilon|,\nabla|u^\varepsilon|\ra-\Delta \eta e^\varepsilon(|u^\varepsilon|)\right)d\vx ds\nn
&+\lim_{h\to0}\lim_{\varepsilon\to 0}\int_{\R} \zeta^h(t-s)\int_{\T}\la\Hess(\eta)\left(\frac{\vj(u^\varepsilon)}{|u^\varepsilon|}-\vj(u^*)\right),\left(\frac{\vj(u^\varepsilon)}{|u^\varepsilon|}-\vj(u^*)\right)\ra d\vx ds\nn&- \lim_{h\to0}\lim_{\varepsilon\to 0}\int_{\R} \zeta^h(t-s)\int_{\T}\Delta \eta\left|\frac{\vj(u^\varepsilon)}{|u^\varepsilon|}-\vj(u^*)\right|^2d\vx ds,
\end{align}
\begin{align}\label{eq:def of I2}
I_2=\lim_{h\to0}\lim_{\varepsilon\to 0}\frac{1}{2}\int_{\R} \zeta^h(t-s)\int_{\T}k_\varepsilon\Delta \eta|u^\varepsilon_t|^2d\vx ds,
\end{align}
\begin{align}
I_3=&\lim_{h\to0}\lim_{\varepsilon\to 0}\int_{\R} \zeta^h(t-s)\int_{\T}2\left(\la\Hess(\eta)\left(\frac{\vj(u^\varepsilon)}{|u^\varepsilon|}-\vj(u^*)\right),\vj(u^*)\ra\right)d\vx ds\nn
&-\lim_{h\to0}\lim_{\varepsilon\to 0}\int_{\R} \zeta^h(t-s)\int_{\T}\Delta\eta\left(\frac{\vj(u^\varepsilon)}{|u^\varepsilon|}-\vj(u^*)\right)\cdot\vj(u^*)d\vx ds.
\end{align}
\eqref{eq:def of j and j/|| NW} immediately implies 
\begin{equation}\label{eq:I_3=0 NW}
I_3=0.
\end{equation}
To estimate $I_1$, combining \eqref{eq:hamiltonian conservation NW} and \eqref{con:limit of integer of j(varphi),energy} , we have

\vspace{-0.35cm}
\begin{align}\label{eq:estimate of Eu W(a0) NW}
    E^\varepsilon(u^\varepsilon(t))=&\int_{\T}h^\varepsilon(u^\varepsilon(\vx,t))d\vx -\int_{\T}\frac{k_\varepsilon}{2}|u^\varepsilon_t|^2d\vx \le W_\varepsilon(\va^0)-\int_{\T}\frac{k_\varepsilon}{2}|u^\varepsilon_t|^2d\vx+o(1).
\end{align}
Taking the derivative of $W(\va(t))+\frac{\pi}{2}|\dot{\va}(t)|^2$ with respect to $t$ and noting \eqref{eq:NWODE}, we have

\vspace{-0.35cm}
\begin{align*}
&\frac{d}{dt}\left(W(\va(t))+\frac{\pi}{2}|\dot{\va}(t)|^2\right)=\sum_{j=1}^{2N}\left(\dot{\va}_j(t)\cdot\nabla_{\va_j}W(\va(t))+\pi\dot{\va}_j(t)\cdot\ddot{\va}_j(t)\right)\\
&\quad=\sum_{j=1}^{2N}\left(\dot{\va}_j(t)\cdot\nabla_{\va_j}W(\va(t))-\dot{\va}_j(t)\cdot\nabla_{\va_j}W(\va(t))\right)=0,
\end{align*}
which immediately implies
\begin{equation}\label{eq:conserve of W NW}
W(\va(t))+\frac{\pi}{2}|\dot{\va}(t)|^2\equiv W(\va(0)).
\end{equation}
Combining \eqref{eq:conserve of W NW}, \eqref{eq:def of We}, \eqref{eq:Lip W} and \eqref{eq:estimate of Eu W(a0) NW}, we have

\vspace{-0.35cm}
\begin{align*}
E^\varepsilon(u^\varepsilon(t))\le &W_\varepsilon(\va(t))+\frac{\pi}{2}|\dot{\va}(t)|^2 -\int_{\T}\frac{k_\varepsilon}{2}|u^\varepsilon_t|^2d\vx +o(1)\\
    \le&W_\varepsilon(\vb(t))+C\zeta(t)+\frac{\pi}{2}|\dot{\va}(t)|^2-\int_{\T}\frac{k_\varepsilon}{2}|u^\varepsilon_t|^2d\vx +o(1),
\end{align*}
which together with \eqref{eq:def of I1}, \eqref{eq:def and gradient of eta}, \eqref{eq:lower bound of L2 ut} and \eqref{eq:distance between ju and ju*} implies

\vspace{-0.35cm}
\begin{align}\label{eq:estimate of I1 NW}
I_1\le&\lim_{h\to0}\lim_{\varepsilon\to0}C\frac{1}{h}\int_{t-h}^{t+h}\int_{B_{r_*}(\vb_j(t))\setminus B_{3r_*/4}(\vb_j(t))}\left(\left|\frac{\vj(u^\varepsilon)}{|u^\varepsilon|}-\vj(u^*) \right|^2+e^\varepsilon(|u^\varepsilon|)\right)d\vx ds\nn
\le&\lim_{h\to0}\lim_{\varepsilon\to0}C\frac{1}{h}\int_{t-h}^{t+h}\left(\frac{\pi}{2}|\dot{\va}(s)|^2-\int_{\T}\frac{k_\varepsilon}{2}|u^\varepsilon_t|^2d\vx +C\zeta(s) \right)ds\nn
\le&\lim_{h\to0}C\frac{1}{h}\int_{t-h}^{t+h}\left(\frac{\pi}{2}|\dot{\va}(s)|^2-\frac{\pi}{2}|\dot{\vb}(s)|^2\right)ds+C\zeta(t)\le C\zeta(t).
\end{align}

Then we estimate $I_2$. Combining \eqref{eq:hamiltonian conservation NW}, \eqref{con:limit of integer of j(varphi),energy}, \eqref{eq:conserve of W NW}, \eqref{eq:refined lower bound of e} and \eqref{eq:def of We}, we have

\vspace{-0.35cm}
\begin{align*}
    &\int_{B_{r_*}(\vb_j(t))\setminus B_{3r_*/4}(\vb_j(t))}k_\varepsilon|u^\varepsilon_t(\vx,s)|^2d\vx\\&\quad \le 2(H^\varepsilon(u^\varepsilon(s))-E^\varepsilon(u^\varepsilon(s))) -\sum_{j=1}^{2N}\int_{R_{3r_*/4}(\vb_j(t))}k_\varepsilon|u^\varepsilon_t(\vx,s)|^2d\vx \\
    &\quad\le2W_\varepsilon(\va(s))+\pi|\dot{\va}(s)|^2-2W_\varepsilon(\vb(s))-\pi|\dot{\vb}(s)|^2\\&\qquad+\pi|\dot{\vb}(s)|^2-\sum_{j=1}^{2N}\int_{R_{3r_*/4}(\vb_j(t))}k_\varepsilon|u^\varepsilon_t(\vx,s)|^2d\vx+o(1)\\
    &\quad\le  C\zeta(s)+\pi|\dot{\vb}(s)|^2-\sum_{j=1}^{2N}\int_{R_{3r_*/4}(\vb_j(t))}k_\varepsilon|u^\varepsilon_t(\vx,s)|^2d\vx+o(1),
\end{align*}
which together with \eqref{eq:def of I2},  \eqref{eq:def and gradient of eta} and \eqref{eq:lower bound of L2 ut} implies 

\vspace{-0.35cm}
\begin{align}\label{eq:estimate of I2 NW}
I_2\le&\limsup_{h\to 0}\limsup_{\varepsilon\to 0}\frac{C}{h}\int_{t-h}^{t+h}\int_{B_{r_*}(\vb_j(t))\setminus B_{3r_*/4}(\vb_j(t))}k_\varepsilon|u^\varepsilon_t(\vx,s)|^2d\vx ds\nn
\le&\limsup_{h\to 0}\limsup_{\varepsilon\to 0}\frac{C}{h}\int_{t-h}^{t+h}\left(C\zeta(s)+\pi|\dot{\vb}(s)|^2-\sum_{j=1}^{2N}\int_{R_{3r_*/4}(\vb_j(t))}k_\varepsilon|u^\varepsilon_t(\vx,s)|^2d\vx  \right)ds\nn 
\le&\limsup_{h\to 0}\frac{C}{h}\int_{t-h}^{t+h}(C\zeta(s)+\pi|\dot{\vb}(s)|^2-\pi|\dot{\vb}(s)|^2)ds
\le C\zeta(t).
\end{align}
Substituting \eqref{eq:splitting of A NW}, \eqref{eq:estimate of I1 NW}, \eqref{eq:estimate of I2 NW} and \eqref{eq:I_3=0 NW} into \eqref{eq:first estimate of zeta' NW}, we obtain 
\[
    \dot{\zeta}(t)\le C\zeta(t),t\in[0,T_*],
\]
which together with \eqref{eq:zeta0=0} implies 
\begin{equation}\label{eq:zetat=0}
    \zeta(t)=0,t\in[0,T_*].
\end{equation}
Combining \eqref{eq:def of zeta} and \eqref{eq:zetat=0}, we obtain $\va(t)=\vb(t)$, which together with \eqref{eq:converge of Ju and eu NW} and that $\va$ satisfies \eqref{eq:NWODE} finishes the proof.
\end{proof}

\section{Numerical simulations of the dynamics of vortex dipole}\label{sec:Numerical}
The momentum $\bvec{Q}(u^\varepsilon_0)$ makes great contribution to the vortex motion of \eqref{eq:NW}. With Theorem \ref{thm:dynamics NW}, we can study the influence of the momentum on vortex motion of \eqref{eq:NW} with $0<\varepsilon<<1$ by simulating \eqref{eq:NWODE} with initial data \eqref{eq:initial ODE} satisfying $\J\vq_0=\lim_{\varepsilon\to 0}\bvec{Q}(u^\varepsilon_0)$.

\begin{figure}[!htp]
\begin{center}
\begin{tabular}{cc}
\includegraphics[height=5.5cm]{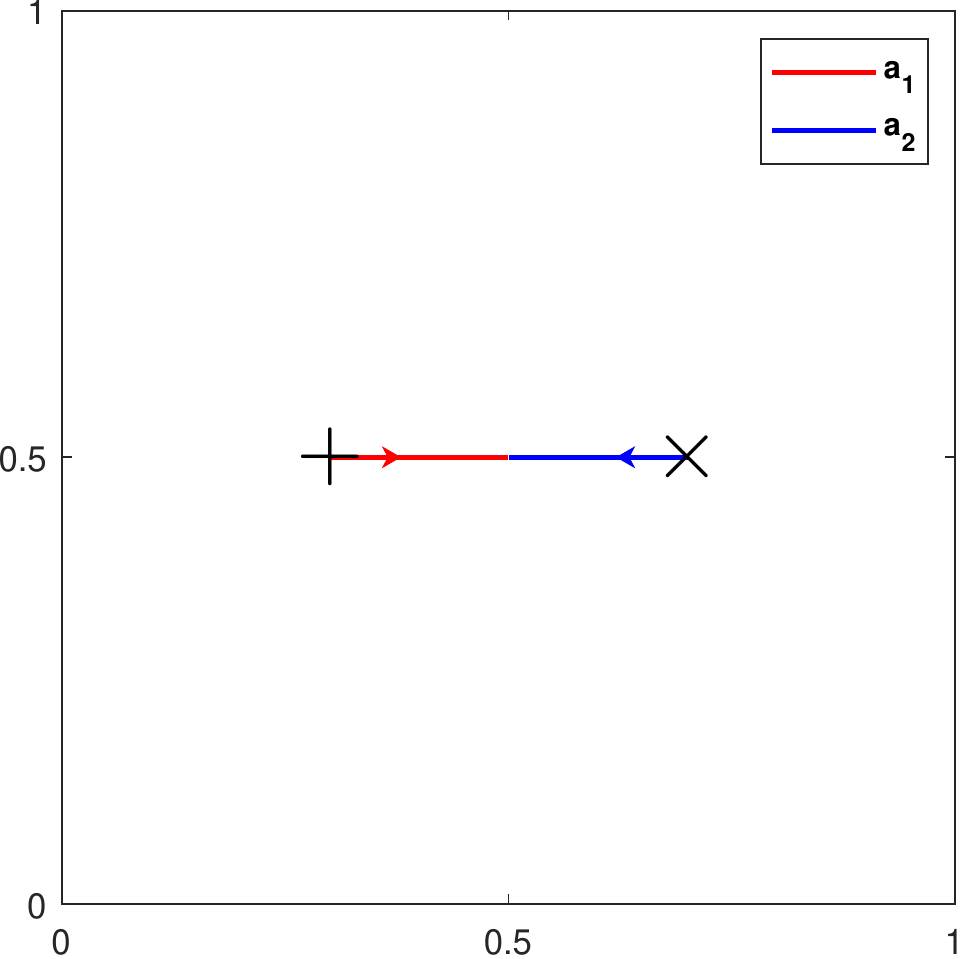}
&\includegraphics[height=5.5cm]{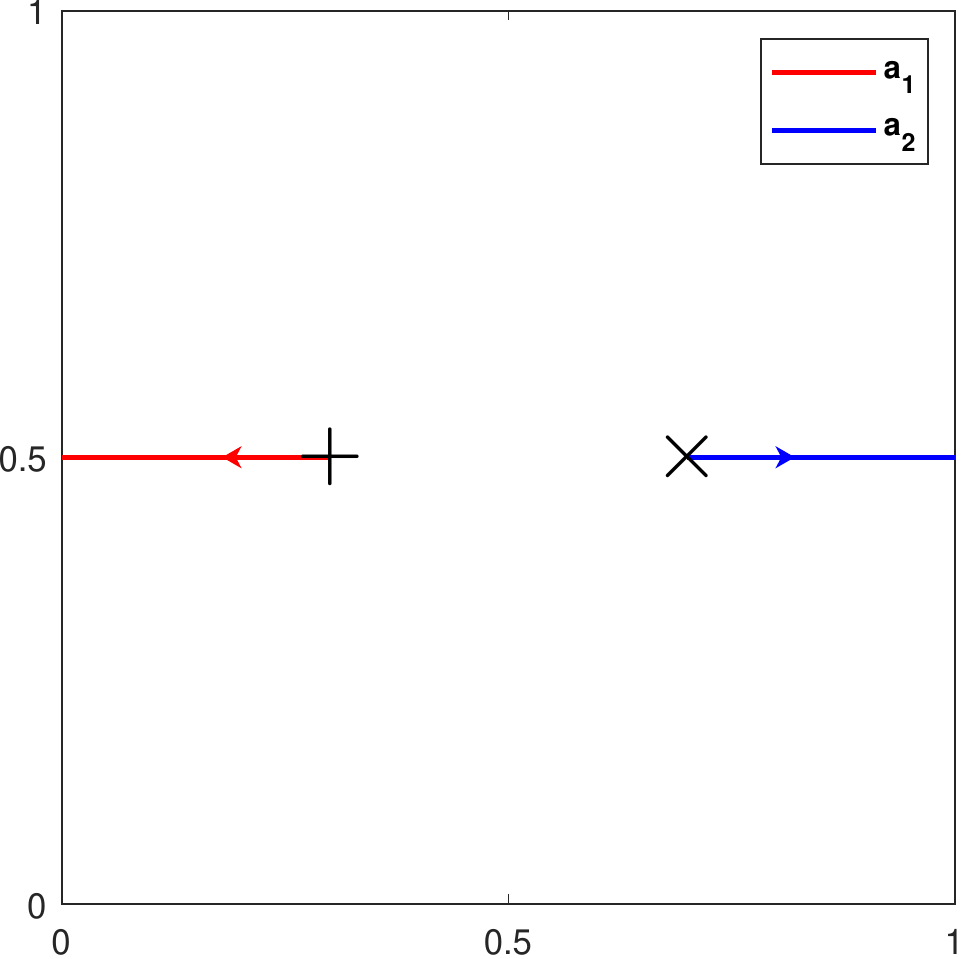}
\end{tabular}
\end{center}
\caption{Trajectories of \eqref{eq:NWODE} with $N=1$ and different initial datum in \eqref{eq:initial ODE} for: (i) $\va_1^0=(0.3,0)^T,\va_2^0=(0.7,0)^T,\vq_0=2\pi(\va_1^0-\va_2^0)$ (left), (ii) $\va_1^0=(0.3,0)^T,\va_2^0=(0.7,0)^T,\vq_0=2\pi(\va_1^0-\va_2^0+(1,0)^T)$ (right). Here and in the below, we use $+$ and $\times$ in the pictures to denote vortices with degrees $+1$ and $-1$, respectively.}\label{fig:paths}
\end{figure}

\begin{figure}[htp!]
\begin{center}
\begin{tabular}{cc}
\includegraphics[height=5.5cm]{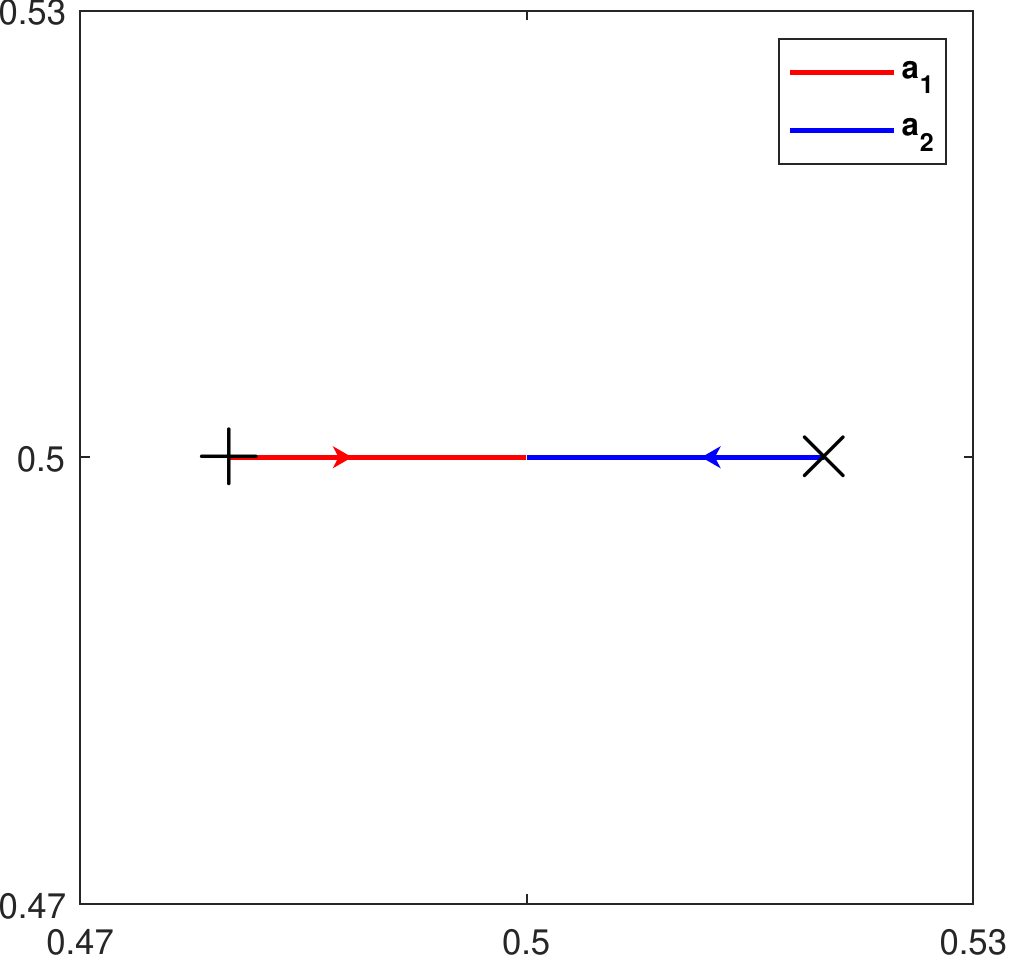}
&\includegraphics[height=5.5cm]{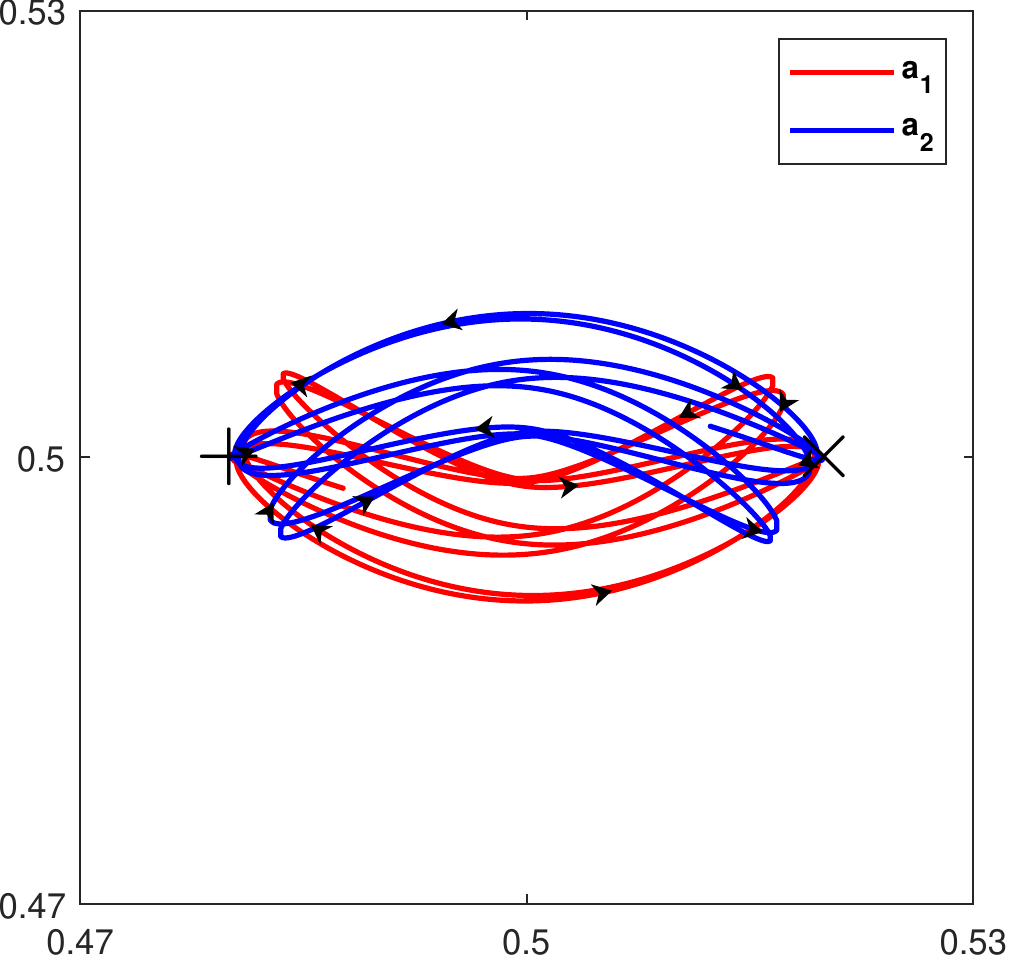}\\
\includegraphics[width=5.5cm]{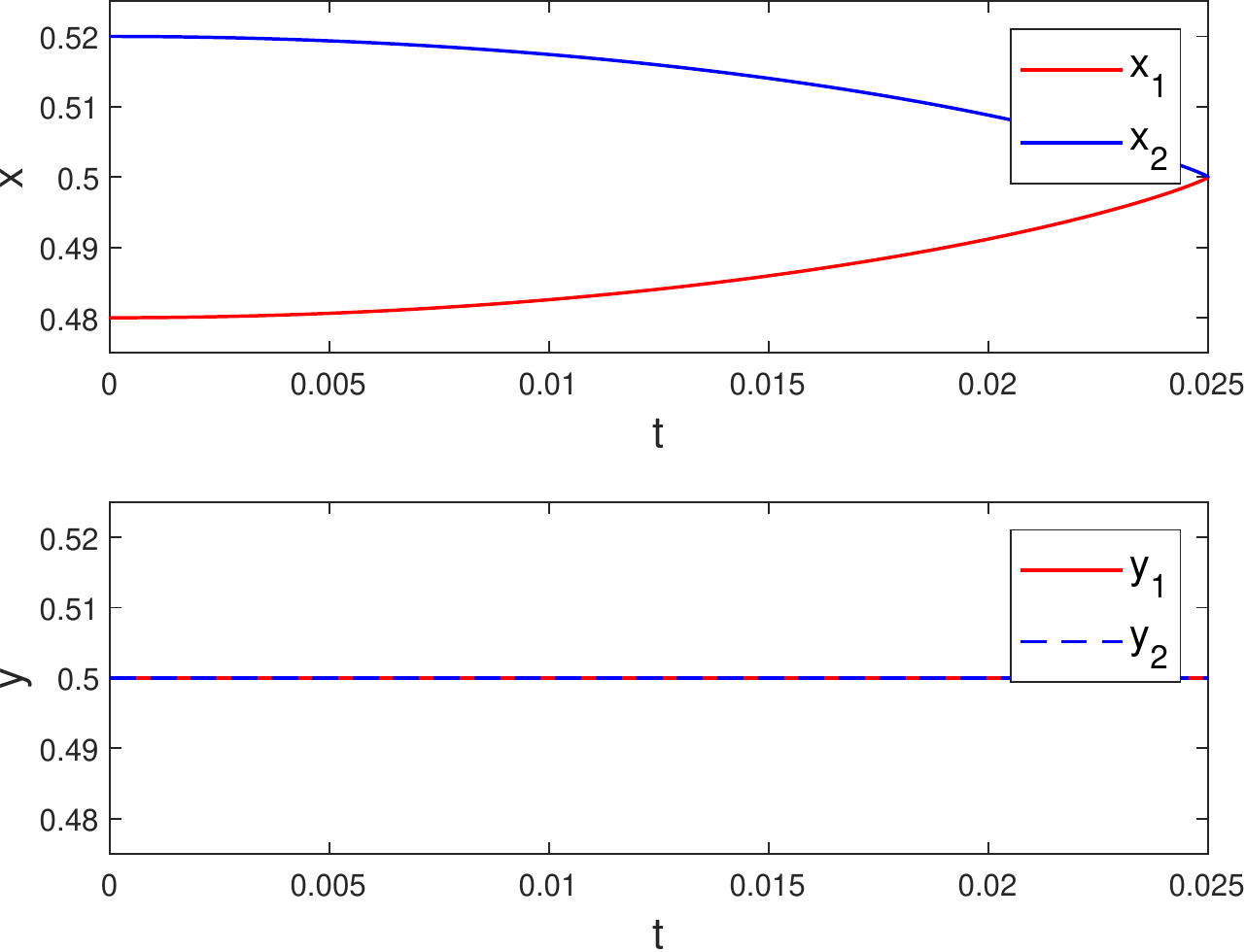}
&\includegraphics[width=5.5cm]{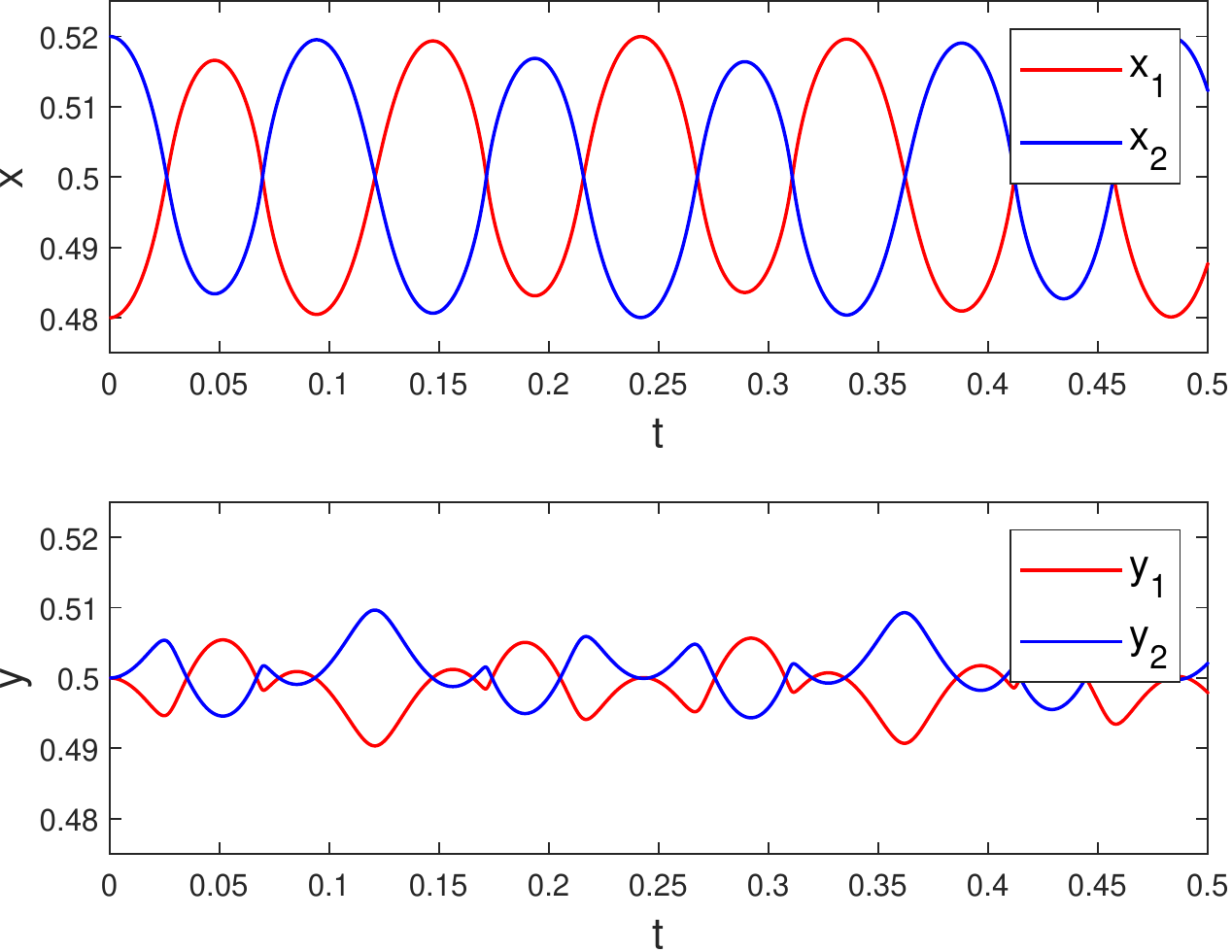}
\end{tabular}
\end{center}
\caption{Trajectories (upper line) and values of  $\va_1=(x_1,y_1)^T,\va_2=(x_2,y_2)^T$ (lower line) of \eqref{eq:NWODE} with $N=1$ and different initial datum in \eqref{eq:initial ODE} for: (i) $\va_1^0=(0.48,0)^T,\va_2^0=(0.52,0)^T,\vq_0=2\pi(\va_1^0-\va_2^0)$ (left), (ii) $\va_1^0=(0.48,0)^T,\va_2^0=(0.52,0)^T,\vq_0=2\pi(\va_1^0-\va_2^0+(0,2)^T)$ (right).}\label{fig:t-x}
\end{figure}

We focus on the case of the vortex dipole, i.e. $N=1$. To illustrate the solution of \eqref{eq:NWODE} with different initial data \eqref{eq:initial ODE}, we solve it numerically by adopting the fourth-order Runge-Kutta method with time step $\Delta t=5\times 10^{-6}$. The results are recorded in Figure \ref{fig:paths} and Figure \ref{fig:t-x}.


\section*{Acknowledgements} The author would like to express his sincere gratitude to Prof. Huaiyu Jian
in Tsinghua University and Prof. Weizhu Bao in National University of Singapore for their guidance and
encouragement.

\end{document}